\documentclass[a4paper,12pt,reqno]{amsart}
\usepackage{latexsym}
\usepackage{amssymb,amsfonts,amsmath,mathrsfs}
\newtheorem{theorem}{Theorem}[section]

\newtheorem{lemma}{Lemma}[section]
\newtheorem{conjecture}{Conjecture}[section]

\author{H. M. Bui, Steven M. Gonek \and Micah B. Milinovich}
\address{Institut f\"ur Mathematik, Universit\"at Z\"urich, Z\"urich CH-8057, Switzerland}
\email{hung.bui@math.uzh.ch}
\address{Department of Mathematics, University of Rochester, Rochester, NY 14627, USA}
\email{gonek@math.rochester.edu} 
\address{Mathematics Department, The University of Mississippi, University, MS 38677, USA}
\email{mbmilino@olemiss.edu}


\title[A hybrid Euler-Hadamard product and moments of $\zeta'(\rho)$]{A hybrid Euler-Hadamard product \\ and moments of $\zeta'(\rho)$}
\begin{document}
\begin{abstract}
Keating and Snaith modeled the Riemann zeta-function $\zeta(s)$ by characteristic polynomials of  random  $N\times N$ unitary matrices, and used this to conjecture the asymptotic main term for the $2k$-th moment of $\zeta(1/2+it)$ when $k>-1/2$. However, an arithmetical factor, widely believed to be part of the leading term coefficient, had to be inserted  in an  \emph{ad hoc} manner. Gonek, Hughes and Keating  later developed a hybrid  formula for $\zeta(s)$ that  combines a truncation of its Euler product  with a  product over its zeros.  Using it, they recovered the moment conjecture of Keating and Snaith in a way that naturally includes the arithmetical factor. Here we use the hybrid formula to recover a conjecture of Hughes, Keating and O'Connell concerning discrete moments of the derivative of the Riemann zeta-function averaged over the zeros of $\zeta(s)$, incorporating the arithmetical factor in a natural way.
\end{abstract}
\maketitle

\section{Introduction}

Let $\zeta(s)$ denote the Riemann zeta-function. In this paper, we study discrete moments of $\zeta'(s)$ in the form
\begin{equation*}
J_k(T)=\frac{1}{N(T)}\sum_{0<\gamma\leq T} \big|\zeta'(\rho) \big|^{2k},
\end{equation*}
where the summation is over the non-trivial zeros $\rho=\beta+i\gamma$ of $\zeta(s)$, and $N(T)$ is the usual zero counting function
\begin{equation*}
N(T)=\sum_{0<\gamma\leq T}1=\frac{T\mathscr{L}}{2\pi}-\frac{T}{2\pi}+O(\mathscr{L}).
\end{equation*}
Here and throughout the paper, we let $\mathscr{L}=\log \frac{T}{2\pi}$, and all sums involving the zeros of $\zeta(s)$ are counted with multiplicity.

The function $J_k(T)$ is defined for all $k\geq0$, and, on the additional assumption that all the zeros are simple, for all $k\in \mathbb{R}$. Trivially, $J_0(T)=1$, but it is still an open problem to rigorously determine the behavior of $J_k(T)$ for any other value of $k$. Gonek~\cite{G3} proved that if the Riemann Hypothesis (RH) is true, then  $J_1(T)\sim \frac{1}{12}\mathscr{L}^{3}$ as $T \to \infty$. Conrey and Snaith~\cite{CS} conjectured the full asymptotic formula for $J_1(T)$ using the $L$-functions Ratios Conjecture, and Milinovich~\cite{M1} proved that  their formula is correct assuming RH.

For   $k$ in general,   Gonek~\cite{G2} and Hejhal~\cite{H} independently conjectured that
\begin{equation}\label{220}
J_k(T)\asymp_k\mathscr{L}^{k(k+2)}
\end{equation}
for fixed $k\in\mathbb{R}$, as $T\rightarrow\infty$. This conjecture is widely believed for non-negative values of $k$, but there is evidence that it is false for $k\leq-3/2$. The case $k=1$ of \eqref{220} holds on RH, of course, by the remarks above, and Ng~\cite{Ng} established the case $k=2$ assuming RH. The conjectured lower bound is known to hold for $k=-1$ under the additional condition that all the zeros of $\zeta(s)$ are simple [\ref{G2}, \ref{MNg2}], and for all $k\in\mathbb{N}$ assuming the generalized Riemann Hypothesis for Dirichlet $L$-functions \cite{MNg}. Moreover, Milinovich~\cite{M} also proved that the upper bound
\begin{equation*}
J_{k}(T)\ll_{k,\varepsilon}\mathscr{L}^{k(k+2)+\varepsilon}
\end{equation*}
holds for all fixed $k\in\mathbb{N}$ and any $\varepsilon>0$ on RH.


The conjecture of Gonek and Hejhal has been refined further using random matrix theory. Let $U$ denote an $N\times N$ unitary matrix  with eigenangles $\theta_n$ $(n=1, 2, \ldots,N)$, and denote its characteristic polynomial by
\begin{equation*}
Z(\theta)=\textrm{det}\big(I-Ue^{-i\theta}\big)=\prod_{n=1}^{N}\big(1-e^{i(\theta_n-\theta)}\big).
\end{equation*}
The random matrix theory model for $J_k(T)$ is
\begin{equation}\label{223}
\int_{U(N)}\frac{1}{N}\sum_{n=1}^{N}\big|Z'(\theta_n)\big|^{2k}d\mu_N,
\end{equation}
where the integral is  over all    $N\times N$ unitary matrices with respect to Haar measure. 
Hughes, Keating and O'Connell [\ref{HKO}] showed that this expression is equal to
\begin{equation}\label{222}
\frac{G^2(k+2)}{G(2k+3)}\frac{G(N)G(N+2k+2)}{NG^2(N+k+1)}\sim\frac{G^2(k+2)}{G(2k+3)}N^{k(k+2)}
\end{equation}
for any fixed $k$ with $\Re(k)>-3/2$, as $N\rightarrow\infty$. Here $G(k)$ is the Barnes $G$-function. Equating the mean densities of the zeros of $\zeta(s)$ and the eigenangles of $U$, that is to set
\begin{equation*}
N\sim\mathscr{L},
\end{equation*}
they were led to the following conjecture.

\begin{conjecture} \label{HKO_conj}\textbf{\textup{(Hughes, Keating and O'Connell)}}
For any fixed $k$ with $\Re(k)>-3/2$, we have
\begin{equation*}
J_k(T)\sim a_k\frac{G^2(k+2)}{G(2k+3)}\mathscr{L}^{k(k+2)}
\end{equation*}
as $T\rightarrow\infty$, where
\begin{equation}\label{221}
a_k=\prod_{p\ \emph{prime}}\bigg(1-\frac{1}{p}\bigg)^{k^2}\sum_{m=0}^{\infty}\bigg(\frac{\Gamma(m+k)}{m!\Gamma(k)}\bigg)^2p^{-m}.
\end{equation}
\end{conjecture}

We note that this agrees with the result  $J_1(T)\sim \frac{1}{12}\mathscr{L}^3$ proved by Gonek~\cite{G3} on RH, and also recovers a conjecture of Gonek \cite{G2,G4} in the case $k=-1$. 
The work of Hughes, Keating and O'Connell is closely related to the work of Keating and Snaith~\cite{KS}, in which they used the characteristic polynomials of large random unitary matrices to model the value distribution of the Riemann zeta-function and study the moments of $\zeta(1/2+it)$. Evaluating the moments of $|Z(\theta)|$ over $U(N)$ with respect to Haar measure and setting $N\sim\mathscr{L}$, they made the following conjecture.

\begin{conjecture}\label{KS_conj} \textbf{\textup{(Keating and Snaith)}}
For any fixed $k$ with $\Re(k)>-1/2$, we have
\begin{equation*}
\frac{1}{T}\int_{0}^{T}\big|\zeta(\tfrac{1}{2}+it)\big|^{2k}\sim a_k\frac{G^2(k+1)}{G(2k+1)}\mathscr{L}^{k^2}
\end{equation*}
as $T\rightarrow\infty$, where $a_k$ is defined as in \eqref{221}.
\end{conjecture}

In both Conjecture \ref{HKO_conj} and Conjecture \ref{KS_conj}, the arithmetical factor $a_k$ was inserted in an \emph{ad hoc} manner  based upon separate number theoretic considerations. This is a  typical drawback of random matrix models of the Riemann zeta-function and other  $L$-functions: they contain no arithmetical information. Moreover, there is no explanation as to why the arithmetical factor $a_k$ is the same in both conjectures; indeed continuous averages of Dirichlet polynomials and averages of Dirichlet polynomials over the zeros of $\zeta(s)$ behave differently.
 
 Gonek, Hughes and Keating~\cite{GHK} developed a new model for $\zeta(s)$ that incorporates the arithmetical information in a natural way.   Their ``hybrid'' model  is based on an approximation of the Riemann zeta-function at a height $t$ on the critical line by a partial Euler product, $P_X(1/2+it)$, multiplied by what is essentially a partial Hadamard product, $Z_X(1/2+it)$, over the non-trivial zeros of $\zeta(s)$ close to $1/2+it$ (see the definitions of $P_X(s)$ and $Z_X(s)$ in the next section). That is,  $\zeta(s)$ is represented as a product over a finite number of primes and   zeros. The moments of  $P_X(s)$ can be calculated rigorously and give rise to the arithmetical factor $a_k$, whereas the moments of the truncated Hadamard 
 product are conjectured using random matrix theory. Under the assumption that the moments  of $\zeta(s)$ split as the product of the moments of $P_X(s)$ and $Z_X(s)$, which can be  proved in certain   cases, they again arrived at Conjecture \ref{KS_conj}. An interesting feature of their approach is that the arithmetic and random matrix theory aspects are treated on an equal footing.  Subsequently, the hybrid Euler-Hadamard product has been extended to various families of 
 $L$-functions \cite{BK1, BK2, D}.

In this paper, we adapt Gonek, Hughes and Keating's model to the problem of estimating $J_k(T)$.  As before, our calculations  suggest that the discrete moments of the derivative of the Riemann zeta-function are asymptotic to the discrete moments of $P_X(s)$ times the discrete moments of the derivative of $Z_X(s)$. Moreover, the model explains why the same arithmetical factor $a_k$ appears in both Conjecture \ref{HKO_conj} and Conjecture \ref{KS_conj}, above.

\section{Hybrid Euler-Hadamard product and the main results}

We begin by stating the hybrid Euler-Hadamard product formula of Gonek, Hughes and Keating (Theorem 1 of \cite{GHK}).

\begin{theorem}\label{thm 1}
Let $X\geq 2$ and $f$ be a non-negative $C^\infty$-function of mass $1$ supported on $[0,1]$. Define 
\begin{equation*}
U(z)=\int_{0}^{1}f(u)E_1\big(z(u\!+\!X\!-\!1)/X\big)du,
\end{equation*}
where $E_1(z)=\int_{z}^{\infty}e^{-u}/u \, du$ is the exponential integral.
Then for $\Re(s)=\sigma\geq0$ we have
\begin{equation}\label{213}
\zeta(s)=P_X(s)Z_X(s)\bigg(1+O_{f,B}\bigg(\frac{X^{B+2}}{\big((|s|\!+\!1)\log X\big)^B}\bigg)
+O_f(X^{-\sigma}\log X)\bigg)
\end{equation}
for any $B>0$, where
\begin{equation*}
P_X(s)=\exp\bigg(\sum_{n\leq X}\frac{\Lambda(n)}{n^s\log n}\bigg),
\end{equation*}
$\Lambda(n)$ is the von Mangoldt function, and
\begin{equation*}
Z_X(s)=\exp\bigg(-\sum_{\rho}U\big((s\!-\!\rho)\log X\big)\bigg).
\end{equation*}
\end{theorem}

As was mentioned in \cite{GHK}, $P_X(s)$ is roughly $\prod_{p\leq X}(1-p^{-s})^{-1}$, and $U(z)$ is roughly $E_1(z)$, which is asymptotic to $-\gamma_0-\log z$ for $|z|$ small, where $\gamma_0$ is Euler's constant. Thus, Theorem \ref{thm 1} says that $\zeta(s)$ looks roughly like
\begin{equation*}
\prod_{p\leq X}\bigg(1-\frac{1}{p^{s}}\bigg)^{-1}\prod_{\substack{\rho\\|s-\rho|\ll1/\log X}}\big((s\!-\!\rho)e^{\gamma_0}\log X\big),
\end{equation*}
which is a hybrid formula in that it combines a partial Euler product and (essentially) a partial Hadamard product.  

We note that from the series expansion of $E_1(z)$, we can interpret $\exp(-U(z))$ to be asymptotic to $Cz$ for some constant $C$ as $|z|\rightarrow0$. Hence both $\zeta(s)$ and $Z_X(s)$ vanish at the zeros of the Riemann zeta-function. Using Cauchy's integral formula in a familiar way, we can differentiate both sides of \eqref{213} and maintain an asymptotic formula. In this way, assuming RH, we obtain that 
\begin{equation}\label{225}
\zeta'(\rho)=P_X(\rho)Z_X'(\rho)\bigg(1+O_{f,B}\bigg(\frac{X^{B+2}}{(|\rho|\log X)^B}\bigg)+O_f(X^{-1/2}\log X)\bigg)
\end{equation}
for every non-trivial zero $\rho$ of $\zeta(s)$ (since the term $P_X'(\rho)Z_X(\rho)$ vanishes).

In Section 4, we evaluate the moments of $P_X(\rho)$ rigorously and establish the following theorem.

\begin{theorem}\label{thm 2}
Assume RH. Let $\varepsilon>0$ and $X$, $T\rightarrow\infty$ with $X=O((\log T)^{2-\varepsilon})$. Then for any $k\in\mathbb{R}$ we have
\begin{displaymath}
\frac{1}{N(T)}\sum_{0<\gamma\leq T}\big|P_{X}(\rho)\big|^{2k}=a_k(e^{\gamma_0}\log X)^{k^2}\big(1+O_k\big((\log X)^{-1}\big)\big).
\end{displaymath}
\end{theorem}

Heuristically, we have
\begin{equation*}
Z_X(s)\approx\prod_{\rho}\big((s\!-\!\rho)e^{\gamma_0}\log X\big).
\end{equation*}
Hence
\begin{equation}\label{224}
Z_X'(\rho)\approx(e^{\gamma_0}\log X)W_X(\tilde{\rho}),
\end{equation}
where $\tilde{\rho}=\rho e^{\gamma_0}\log X$, and
\begin{equation*}
W_X(\tilde{\rho})=\prod_{\tilde{\rho}'\ne\tilde{\rho}}\big(\tilde{\rho}-\tilde{\rho}'\big).
\end{equation*}
As in the random matrix model \eqref{223} for $\zeta'(\rho)$ of Hughes, Keating and O'Connell, we  model the $2k$-th moment of $W_X(\tilde{\rho})$ by
\begin{equation*}
\int_{U(N)}\frac{1}{N}\sum_{n=1}^{N}\big|Z'(\theta_n)\big|^{2k}d\mu_N.
\end{equation*}
Here, however,  the average gap  between consecutive  $\tilde{\rho}$'s is $2\pi e^{\gamma_0}\log X/\mathscr{L}$. Therefore, equating the mean density 
of $\tilde{\rho}$ and the density of the eigenangles corresponds to the identification $N\sim\mathscr{L}/e^{\gamma_0}\log X$. Combining   \eqref{222} 
and \eqref{224} leads to the following conjecture.

\begin{conjecture}\label{con 3}
Let $\varepsilon>0$ and $X$, $T\rightarrow\infty$ with $X=O((\log T)^{2-\varepsilon})$. Then for any $k>-3/2$ we have
\begin{displaymath}
\frac{1}{N(T)}\sum_{0<\gamma\leq T}\big|Z_{X}'(\rho)\big|^{2k}\sim\frac{G^2(k+2)}{G(2k+3)}(e^{\gamma_0}\log X)^{2k}\bigg(\frac{\mathscr{L}}{e^{\gamma_0}\log X}\bigg)^{k(k+2)}.
\end{displaymath}
\end{conjecture}

In Section 5  we shall prove the case $k=1$ of Conjecture~\ref{con 3}, assuming RH. Since, by \eqref{225}, 
\[
\zeta'(\rho) P_X(\rho)^{-1} = Z_X'(\rho) \big( 1+o(1) \big),
\]
when $\Im(\rho)=\gamma$ is large and $X=O((\log \gamma)^{2-\varepsilon})$, this amounts to proving the following result.  
\begin{theorem}\label{2nd_moment}
Assume RH. Let $\varepsilon>0$ and $X$, $T\rightarrow\infty$ with $X=O((\log T)^{2-\varepsilon})$. Then we have
\begin{displaymath}
\frac{1}{N(T)}\sum_{0<\gamma\leq T}\big|\zeta'(\rho)P_X(\rho)^{-1}\big|^{2}\sim\frac{1}{12}\frac{\mathscr{L}^3}{e^{\gamma_0}\log X}.
\end{displaymath}
\end{theorem}

In Section 6 we shall use the $L$-functions Ratios Conjectures to heuristically derive the asymptotic formula
\[
\frac{1}{N(T)}\sum_{0<\gamma\leq T}\big|\zeta'(\rho)P_{X}(\rho)^{-1}\big|^4\sim\frac{1}{8640}\frac{\mathscr{L}^8}{(e^{\gamma_0}\log X)^{4}},
\]
and thus, as $1/8640=G^2(4)/G(7)$, provide additional evidence for Conjecture~\ref{con 3} in the case $k=2$. 

Our proof of Theorem \ref{2nd_moment} involves replacing $P_X(\rho)^{-1}$ by a short Dirichlet polynomial and then using the method of Conrey, Ghosh and Gonek \cite{CGG1} to estimate the resulting mean-value. However, unlike the proof in \cite{CGG1}, we do not need to assume the generalized Lindel\"{o}f hypothesis (GLH) for Dirichlet $L$-functions.  We circumvent the assumption of GLH by incorporating ideas of Bui and Heath-Brown \cite{BH-B}, who have recently proved the results in \cite{CGG1} assuming only RH.

Our results for the cases $k=1$ and  $k=2$ suggest that at least  when $X$ is not too large relative to $T$, the $2k$-th discrete moment of $\zeta'(\rho)$ is asymptotic to the 
product of the discrete moments of $P_X(\rho)$ and $Z_X'(\rho)$. We believe that  this is true in general, and we make the following conjecture.

\begin{conjecture}\label{con 4}
Let $\varepsilon>0$ and $X$, $T\rightarrow\infty$ with $X=O((\log T)^{2-\varepsilon})$. Then for any $k>-3/2$ we have
\begin{displaymath}
\frac{1}{N(T)}\sum_{0<\gamma\leq T}\big|\zeta'(\rho)\big|^{2k}\sim\bigg(\frac{1}{N(T)}\sum_{0<\gamma\leq T}\big|P_{X}(\rho)\big|^{2k}\bigg)\bigg(\frac{1}{N(T)}\sum_{0<\gamma\leq T}\big|Z_{X}'(\rho)\big|^{2k}\bigg).
\end{displaymath}
\end{conjecture}

By combining Theorem~\ref{thm 2}, Conjecture \ref{con 3}, and Conjecture \ref{con 4}, we recover the conjecture of Hughes, Keating and O'Connell for real values of $k$ satisfying $k > -3/2$, and incorporate the arithmetical factor $a_k$ in a natural way.

\section{Lemmas}\label{sec: lemmas}

In order to prove Theorem \ref{thm 2}, we require the following version of the Landau-Gonek explicit formula~\cite{G}.
\begin{lemma}\label{lem 1}
Let $x,T>1$. Then we have
\begin{equation*}
\begin{split}
\sum_{0<\gamma\leq T}x^{\rho}&=-\frac{T}{2\pi}\Lambda(x)+O(x\log(2xT)\log\log(3x))\\
&\quad\quad +O\bigg(\log x\min\bigg\{T,\frac{x}{\langle x\rangle}\bigg\}\bigg)+O\bigg(\log(2T)\min\bigg\{T,\frac{1}{\log x}\bigg\}\bigg),
\end{split}
\end{equation*}
where $\langle x\rangle$ denotes the distance from $x$ to the nearest prime power other than $x$ itself, and $\Lambda(x)$ is the generalized von Mangoldt function; that is, $\Lambda(x)=\log p$ if $x=p^k$ for a prime $p$ and natural number $k$, and $\Lambda(x)=0$ otherwise.
\end{lemma}

The next two lemmas are in \cite{CGG1} (see Lemma 2 and Lemma 3). 

\begin{lemma}\label{lem 2} 
Suppose that $A(s)=\sum_{m=1}^{\infty}a(m)m^{-s}$, where $a(m)\ll_\varepsilon m^\varepsilon$, and $B(s)=\sum_{n\leq y}b(n)n^{-s}$, where $b(n)\ll_\varepsilon n^\varepsilon$. Then we have
\begin{equation*}
\frac{1}{2\pi i}\int_{c+i}^{c+iT}\chi(1-s)A(s)B(1-s)ds=\sum_{n\leq y}\frac{b(n)}{n}\sum_{m\leq nT/2\pi}a(m)e(-m/n) +O_\varepsilon(yT^{1/2+\varepsilon}),
\end{equation*} 
where $c=1+\mathscr{L}^{-1}$. 
\end{lemma}

\begin{lemma}\label{lem 3}
Suppose that $\alpha=\alpha_1*\alpha_2$. Then we have
\begin{equation*}
\alpha(lm)=\sum_{\substack{l=l_1l_2\\m=m_1m_2\\(m_2,l_1)=1}}\alpha_{1}(l_1m_1)\alpha_{2}(l_2m_2).
\end{equation*}
\end{lemma}

\section{Proof of Theorem \ref{thm 2}}\label{sec: prf thm 2}

Since Theorem \ref{thm 2} holds when $k=0$, we assume throughout this section that $k$ is a nonzero real number. We begin by  approximating  $P_{X}(s)^k$ by a truncated Dirichlet series.  
Write
\begin{equation}\label{2}
P_{X}(s)^k=\sum_{n=1}^{\infty}\frac{\alpha_{k}(n)}{n^s}.
\end{equation}
From  the definition of  $P_{X}(s)$, we see that $\alpha_{k}(n)$ is multiplicative and real valued. Also, if we let 
\begin{displaymath}
S(X)=\{n\in\mathbb{N}:p|n\Rightarrow p\leq X\},
\end{displaymath}
the set of $X$-smooth numbers, then  $\alpha_{k}(n)=0$ if $n\notin S(X)$. 
In \cite{GHK} it is shown that  $|\alpha_{k}(n)|\leq d_{|k|}(n)$, and that $\alpha_{k}(n)=d_{k}(n)$ 
if $n\in S(\sqrt{X})$ or if $n$ is a prime $p\leq X$, where the arithmetic function $d_k(n)$ is defined in terms of the Dirichlet series 
\[
\zeta(s)^k=\sum_{n=1}^{\infty} \frac{d_k(n)}{n^{s}}
\] 
for $\Re(s)>1$ and any real number $k$.  In \cite{GHK} it is also shown (see  page 518) that
\begin{equation}\label{22}
P_{X}(s)^k=\sum_{\substack{n\in S(X)\\n\leq T^{\vartheta}}}
\frac{\alpha_{k}(n)}{n^s}+O_{k,\varepsilon}(T^{-\varepsilon\vartheta/2}) 
\end{equation}
for any $\varepsilon, \vartheta>0$, where   $\vartheta$ will be chosen later.
Using elementary inequalities, we see that
\begin{equation}\label{215}
\begin{split}
\bigg|\bigg(\sum_{0<\gamma\leq T}|P_{X}(\rho)|^{2k}\bigg)^{1/2}&-\bigg(\sum_{0<\gamma\leq T}\bigg|\sum_{\substack{n\in S(X)\\n\leq T^{\vartheta}}}\frac{\alpha_{k}(n)}{n^\rho}\bigg|^{2}\bigg)^{1/2}\bigg|
\\
&\ll_{k,\varepsilon}\bigg(\sum_{0<\gamma\leq T}T^{-\varepsilon\vartheta}\bigg)^{1/2}\ll_{k,\varepsilon}T^{1/2-\varepsilon\vartheta/3}.
\end{split}
\end{equation}
Thus, in order to establish Theorem \ref{thm 2}, it suffices to estimate the second moment of the truncated Dirichlet series. 

Assuming RH, $1-\rho = \overline{\rho}$ for any non-trivial zero $\rho$ of $\zeta(s)$. Therefore
\begin{equation*}
\begin{split}
\sum_{0<\gamma\leq T}\bigg|\sum_{\substack{n\in S(X)\\n\leq T^{\vartheta}}}\frac{\alpha_{k}(n)}{n^\rho}\bigg|^{2} &= \sum_{\substack{mn\in S(X)\\m,n\leq T^\vartheta}}\frac{\alpha_{k}(m)\alpha_{k}(n)}{n}\sum_{0<\gamma\leq T}\bigg(\frac{m}{n}\bigg)^{-\rho} 
\\
&= M+E_1+E_2,
\end{split}
\end{equation*}
say, where $M$, $E_1$, and $E_2$ are the sums representing the contributions from the terms $m=n$, $m<n$, and $m>n$, respectively.
Since $1-\rho = \overline{\rho}$, we see that $E_2=\overline{E_1}$. Thus, it suffices to estimate $E_1$ and $M$. 
From Lemma~\ref{lem 1}, we deduce that $E_1$ equals
\begin{equation*}
\begin{split}
&-\frac{T}{2\pi}\sum_{\substack{mn\in S(X)\\m<n\leq T^\vartheta}}\frac{\alpha_{k}(m)\alpha_{k}(n)}{n}\Lambda\bigg(\frac{n}{m}\bigg) +O\bigg(\mathscr{L}\log\mathscr{L}\sum_{m<n\leq T^\vartheta}\frac{d_{|k|}(m)d_{|k|}(n)}{m}\bigg)\\
&\quad \quad + O\bigg(\mathscr{L}\sum_{m<n\leq T^\vartheta}\frac{d_{|k|}(m)d_{|k|}(n)}{m\langle n/m\rangle}\bigg)+O\bigg(\mathscr{L}\sum_{m<n\leq T^\vartheta}\frac{d_{|k|}(m)d_{|k|}(n)}{n\log n/m}\bigg).
\end{split}
\end{equation*}
We denote these four terms by $E_{11}$, $E_{12}$, $E_{13}$, and $E_{14}$, respectively. Now 
\begin{equation*}
\begin{split}
E_{11}&\ll T\sum_{mn\in S(X)}\frac{d_{|k|}(m)d_{|k|}(n)}{n}\Lambda\bigg(\frac{n}{m}\bigg)
\\
&\ll T\sum_{p\leq X}\sum_{r\geq1}\frac{\log p}{p^r}\sum_{m\in S(X)}\frac{d_{|k|}(m)d_{|k|}(mp^r)}{m}\\
&\ll T\sum_{p\leq X}\sum_{r\geq1}\frac{d_{|k|}(p^r)\log p}{p^r}\sum_{m\in S(X)}\frac{d_{|k|}(m)^{2}}{m}.
\end{split}
\end{equation*}
Since the innermost sum over $m$ is $\ll\prod_{p\leq X}(1-1/p)^{-k^2}\ll_k(\log X)^{k^2}$, it follows that
\begin{displaymath}
E_{11}\ll_k T(\log X)^{k^2}\sum_{p\leq X}\frac{\log p}{p}\ll_k T(\log X)^{k^2+1}.
\end{displaymath}
Trivially we have that
\begin{displaymath}
E_{12}\ll_{k,\varepsilon} T^{\vartheta+\varepsilon}
\end{displaymath}
for any $\varepsilon >0$. To estimate $E_{13}$, we write $n=um+v$ where $|v/m|\leq1/2$. We observe that $\langle n/m\rangle=|v/m|$ if $u$ is a prime power and $v\ne0$, otherwise $\langle n/m\rangle\geq1/2$. Hence
\begin{displaymath}
E_{13}\ll_{k,\varepsilon} T^{\varepsilon}\bigg(\sum_{um\ll T^\vartheta}\sum_{1\leq v\leq m/2}\frac{d_{|k|}(m)}{v}+\sum_{m,n\leq T^\vartheta}\frac{d_{|k|}(m)d_{|k|}(n)}{m}\bigg)\ll_{k,\varepsilon} T^{\vartheta+\varepsilon}.
\end{displaymath}
For $E_{14}$, we note that $\log \frac{n}{m} \geq\log \frac{n}{n-1} \gg 1/n$. Therefore 
\begin{equation*}
E_{14}\ll_{\varepsilon} T^{\varepsilon}\sum_{m,n\leq T^\vartheta}d_{|k|}(m)d_{|k|}(n)\ll_{k,\varepsilon} T^{2\vartheta+\varepsilon}. 
\end{equation*}
Combining the above estimates, we have shown that
\begin{equation}\label{4}
E_1\!+\!E_2 \ll_{k,\varepsilon} T(\log X)^{k^2+1}+T^{2\vartheta+\varepsilon}.
\end{equation}

For the evaluation of $M$,  we appeal to  Lemma~\ref{lem 2} of \cite{GHK} and its proof, and get
\begin{equation}\label{3}
\begin{split}
M&=N(T)\sum_{\substack{n\in S(X)\\n\leq T^\vartheta}}\frac{\alpha_{k}(n)^{2}}{n}
\\
&=N(T)a_k(e^{\gamma_0}\log X)^{k^2}\big(1+O_{k}\big((\log X)^{-1}\big)\big).
\end{split}
\end{equation}
Theorem \ref{thm 2} now follows from \eqref{215}, \eqref{4}, and \eqref{3} by choosing any $\vartheta<1/2$. \\

\noindent{\sc Remark.} The above proof illustrates why the arithmetical factor $a_k$ is the same in both Conjecture \ref{HKO_conj} and Conjecture \ref{KS_conj}, and this arises from a combination of two different phenomena. First of all, while $\zeta'(s)$ is approximated by $P_X'(s)Z_X(s) + P_X(s)Z_X'(s) $, as we noted above $\zeta'(\rho)$ is approximated by $P_X(\rho)Z_X'(\rho) $. Consequently, the arithmetical factor $a_k$ arises solely from moments of the truncated Euler product $P_X(s)$, and not from the moments of its derivative $P_X'(s)$. Moreover, as is the case with continuous moments of $P_X(s)$, there is no off-diagonal contribution to the main term of these moments. For a ``typical" Dirichlet polynomial we expect an additional main term contribution from the sum corresponding to $E_{11}$ in the above proof. However, in the present case, the arithmetic nature of the coefficients $\alpha_k(n)$ (i.e. supported on $X$-smooth numbers with $X=O((\log T)^{2-\varepsilon})$) implies that the term $E_{11}$ contributes an amount which is an error term.


\section{Proof of Theorem \ref{2nd_moment}}

\subsection{Initial setup}

Using the expression in \eqref{22} with $k=-1$,  we have
\begin{equation}\label{201}
\sum_{0<\gamma\leq T}\big|\zeta'(\rho)P_{X}(\rho)^{-1}\big|^2
=\sum_{\substack{mn\in S(X)\\m,n\leq T^{\vartheta}}}\frac{\alpha_{-1}(m)\alpha_{-1}(n)}{\sqrt{mn}}I(m,n)+O_{\varepsilon}\big(T^{1-\varepsilon\vartheta/3}\big),
\end{equation}
where
\begin{equation*}
I(m,n)=\sum_{0<\gamma\leq T}\big|\zeta'(\rho)\big|^2\bigg(\frac{m}{n}\bigg)^{-i\gamma}.
\end{equation*}
Throughout the proof of Theorem \ref{2nd_moment}, we shall repeatedly use the estimate $|\alpha_{-1}(n)| \le d(n)$, where $d(n)$ is the divisor function.
 
We differentiate both sides of the functional equation
\begin{equation*}
\zeta(s)=\chi(s)\zeta(1-s) 
\end{equation*}
to obtain
\begin{equation}\label{301}
\zeta'(s)=-\chi(s)\bigg(\zeta'(1-s)-\frac{\chi'}{\chi}(s)\zeta(1-s)\bigg).
\end{equation}
It follows that
 $\zeta'(1-\rho)=-\chi(1-\rho)\zeta'(\rho)$. 
Thus, assuming RH and using Cauchy's theorem, we get
\begin{eqnarray*}
I(m,n)&=&-\sum_{0<\gamma\leq T}\chi(1-\rho)\zeta'(\rho)^2\bigg(\frac{m}{n}\bigg)^{-i\gamma}\\
&=&-\frac{1}{2\pi i}\int_{\mathscr{C}}\chi(1-s)\frac{\zeta'}{\zeta}(s)\zeta'(s)^2\bigg(\frac{m}{n}\bigg)^{-s+1/2}ds,
\end{eqnarray*}
where $\mathscr{C}$ is the positively oriented rectangle with vertices at 
$1-c+i$, $c+i$, $c+iT$ and $1-c+iT$. Here $c=1+\mathscr{L}^{-1}$ and 
$T$ is chosen so that the distance from $T$ to the nearest ordinate 
of a zero is $\gg \mathscr{L}^{-1}$. 

By standard estimates, for $s$ on $\mathscr{C}$ we have $\zeta'(s)/\zeta(s)\ll \mathscr{L}^2$, 
$\zeta'(s)\ll T^{(1-\sigma)/2}\mathscr{L}$, and $\chi(1-s)\ll T^{\sigma-1/2}$. 
Hence, the contribution from the horizontal segments of $\mathscr{C}$ is 
$$\ll_\varepsilon (m+n)(mn)^{-1/2}T^{1/2+\varepsilon} .$$

We denote the contributions from the right-hand and left-hand edges of $\mathscr{C}$ by $I_R(m,n)$ and $I_L(m,n)$, respectively. Thus,
\begin{equation}\label{100}
I_R(m,n)=-\frac{1}{2\pi i}\int_{c+i}^{c+iT}\chi(1-s)\frac{\zeta'}{\zeta}(s)\zeta'(s)^2\bigg(\frac{m}{n}\bigg)^{-s+1/2}ds,
\end{equation}
and $I_L(m,n)$ is the same except that the  integral is from $1-c+iT$ to $1-c+i$. Logarithmically differentiating the functional equation, we have
\begin{equation}\label{log_deriv}
\frac{\zeta'}{\zeta}(1-s)=\frac{\chi'}{\chi}(1-s)-\frac{\zeta'}{\zeta}(s).
\end{equation}
Using \eqref{301} twice and substituting $1-s$ for $s$, we see that  
\begin{eqnarray*}
I_L(m,n)&=&-\frac{1}{2\pi i}\int_{c-i}^{c-iT}\chi(1-s)\bigg(\frac{\chi'}{\chi}(1-s)-\frac{\zeta'}{\zeta}(s)\bigg)
\\&&  \qquad \qquad\qquad\qquad \times \bigg(\zeta'(s)-\frac{\chi'}{\chi}(1-s)\zeta(s)\bigg)^2\bigg(\frac{m}{n}\bigg)^{s-1/2}ds\nonumber\\
&=&\overline{I_R(n,m)}+\overline{I'(m,n)}+\overline{I''(m,n)},
\end{eqnarray*}
where
\begin{equation*}
I'(m,n)=\frac{1}{2\pi i}\int_{c+i}^{c+iT}\frac{\chi'}{\chi}(1-s)^3\zeta(s)\zeta(1-s)\bigg(\frac{m}{n}\bigg)^{s-1/2}ds
\end{equation*}
and
\begin{equation*}
I''(m,n)=\frac{-3}{2\pi i}\int_{c+i}^{c+iT}\frac{\chi'}{\chi}(1-s)\zeta'(s)\zeta'(1-s)\bigg(\frac{m}{n}\bigg)^{s-1/2}ds.
\end{equation*}
Thus,
\begin{equation*}
\begin{split}
I(m,n)= I_R(m,n)+\overline{I_R(n,m)}+\overline{I'(m,n)}+\overline{I''(m,n)}+O_\varepsilon\big((m+n)(mn)^{-1/2}T^{1/2+\varepsilon}\big).
\end{split}
\end{equation*}
 We shall write the sum on the right-hand side of \eqref{201} as 
 \begin{equation}\label{J sum}
 \sum_{\substack{mn\in S(X)\\m,n\leq T^{\vartheta}}}\frac{\alpha_{-1}(m)\alpha_{-1}(n)}{\sqrt{mn}}I(m,n)
 =J_1+J_2+J_3+J_4+J_5
 \end{equation}
corresponding to this decomposition of $I(m,n)$.

\subsection{The evaluation of $J_3$, $J_4$ and $J_5$}

The term $J_5$ is easy to handle since
\begin{equation}\label{202}
J_5\ll_\varepsilon T^{1/2+\varepsilon}\sum_{m,n\leq T^\vartheta}\frac{d(m)d(n)(m+n)}{mn}\ll_\varepsilon T^{1/2+\vartheta+\varepsilon}.
\end{equation}
To estimate $J_3$ and $J_4$, we move the line of integration in both $I'(m,n)$ and $I''(m,n)$ to the $\tfrac{1}{2}$-line. As in \eqref{202}, this produces an error of size $O_\varepsilon(T^{1/2+\vartheta+\varepsilon})$. Therefore
\begin{equation}\label{206}
J_3=\frac{1}{2\pi}\int_{1}^{T}\frac{\chi'}{\chi}\big(\tfrac{1}{2}+it\big)^3\big|\zeta(\tfrac{1}{2}+it)\big|^2\bigg|\sum_{\substack{n\in S(X)\\n\leq T^\vartheta}}\frac{\alpha_{-1}(n)}{n^{1/2+it}}\bigg|^2dt+O_\varepsilon(T^{1/2+\vartheta+\varepsilon})
\end{equation}
and
\begin{equation}\label{207}
J_4=-\frac{3}{2\pi}\int_{1}^{T}\frac{\chi'}{\chi}\big(\tfrac{1}{2}+it\big)\big|\zeta'(\tfrac{1}{2}+it)\big|^2\bigg|\sum_{\substack{n\in S(X)\\n\leq T^\vartheta}}\frac{\alpha_{-1}(n)}{n^{1/2+it}}\bigg|^2dt+O_\varepsilon(T^{1/2+\vartheta+\varepsilon}).
\end{equation}

Let
\begin{equation}\label{302}
J_3'=\int_{1}^{T}\big|\zeta(\tfrac{1}{2}+it)\big|^2\bigg|
\sum_{\substack{n\in S(X)\\n\leq T^\vartheta}}\frac{\alpha_{-1}(n)}{n^{1/2+it}}\bigg|^2dt
\end{equation}
and
\begin{equation}\label{303}
J_4'=\int_{1}^{T}\big|\zeta'(\tfrac{1}{2}+it)\big|^2\bigg|
\sum_{\substack{n\in S(X)\\n\leq T^\vartheta}}
\frac{\alpha_{-1}(n)}{n^{1/2+it}}\bigg|^2dt.
\end{equation}
If $\vartheta < \frac{1}{2}$, then the integral in \eqref{302} is of the form evaluated in \cite{BCH}, while the integral in \eqref{303} is almost of this form, but not quite. However,
with obvious changes to the argument in \cite{BCH} that we will not carry out here, one may show that 
\begin{eqnarray}\label{203}
&&\int_{1}^{T}\zeta(\tfrac{1}{2}+it+\alpha)\zeta(\tfrac{1}{2}-it+\beta)\bigg|\sum_{n\leq N}\frac{a(n)}{n^{1/2+it}}\bigg|^2dt\nonumber \\
  &&\qquad =\sum_{m,n\leq N}\frac{a(m)\overline{a(n)}(m,n)^{1+\alpha+\beta}}{mn}\nonumber \\
  &&\qquad\qquad\qquad\times \int_{1}^{T}\bigg(m^{-\beta}n^{-\alpha}\zeta(1\!+\!\alpha\!+\!\beta) 
  +\bigg(\frac{t(m,n)^2}{2\pi}\bigg)^{-\alpha-\beta}m^{\alpha}n^\beta\zeta(1\!-\!\alpha\!-\!\beta)\bigg)dt \nonumber \\
&&\qquad\qquad +O_B(T\mathscr{L}^{-B})+O_\varepsilon(N^2 T^{\varepsilon}), 
\end{eqnarray}
uniformly for $\alpha,\beta\ll \mathscr{L}^{-1}$ and for any $B>0$. We use \eqref{203} to estimate both $J_3'$ and $J_4'$.
Applying it first  to \eqref{302}, we find that
\begin{eqnarray}\label{204}
&&\!\!\!\!\!\!\!\!\!\!\!\!J_3'=T\sum_{\substack{mn\in S(X)\\m,n\leq T^\vartheta}}\frac{\alpha_{-1}(m)\alpha_{-1}(n)(m,n)}{mn}\bigg(\log\frac{T(m,n)^2}{2\pi mn}+2\gamma_0-1\bigg)
\nonumber\\
&&\qquad+O_B(T\mathscr{L}^{-B})+O_\varepsilon(T^{2\vartheta+\varepsilon})
\nonumber\\
&&\!\!\!\!\!=T\mathscr{L}\sum_{\substack{mn\in S(X)\\m,n\leq T^\vartheta}}\frac{\alpha_{-1}(m)\alpha_{-1}(n)(m,n)}{mn}+O\bigg(T\sum_{lmn\in S(X)}\frac{d(lm)d(ln)\log mn}{lmn}\bigg)
\nonumber\\
&&\qquad+O_B(T\mathscr{L}^{-B})+O_\varepsilon(T^{2\vartheta+\varepsilon}).
\end{eqnarray}
The double sum in the main term of \eqref{204} has been evaluated by Gonek, Hughes and Keating (see equations (34)--(38) in \cite{GHK}). The analysis in \cite{GHK} implies that 
\begin{equation*}
T\mathscr{L}\sum_{\substack{mn\in S(X)\\m,n\leq T^\vartheta}}\frac{\alpha_{-1}(m)\alpha_{-1}(n)(m,n)}{mn}=\frac{T\mathscr{L}}{e^{\gamma_0}\log X}\big(1+O\big((\log X)^{-1}\big)\big).
\end{equation*}
The sum in the first big-$O$ term of \eqref{204} is
\begin{eqnarray*}
\sum_{lmn\in S(X)}\frac{d(lm)d(ln)\log mn}{lmn}\ll\sum_{l\in S(X)}\frac{d(l)^2}{l}\bigg(\sum_{n\in S(X)}\frac{d(n)\log n}{n}\bigg)^2.
\end{eqnarray*}
Writing
\begin{equation*}
f(\sigma)=\sum_{n\in S(X)}\frac{d(n)}{n^\sigma}=\prod_{p\leq X}\bigg(1-\frac{1}{p^\sigma}\bigg)^{-2},
\end{equation*}
we see that
\begin{equation}\label{212}
\sum_{n\in S(X)}\frac{d(n)\log n}{n}=-f'(1)=2f(1)\sum_{p\leq X}\frac{\log p}{p-1}\ll (\log X)^3.
\end{equation}
Hence the first big-$O$ term in \eqref{204} is $\ll (\log X)^{10}$. Thus, we have shown that
\begin{equation}\label{208}
J_3'=\frac{T\mathscr{L}}{e^{\gamma_0}\log X}\big(1+O\big((\log X)^{-1}\big)\big)+O_\varepsilon(T^{2\vartheta+\varepsilon}).
\end{equation}
Similarly, applying \eqref{203} to \eqref{303}, we obtain
\begin{equation}\label{209}
\begin{split}
J_4'&=\frac{T\mathscr{L}^3}{3}\sum_{\substack{mn\in S(X)\\m,n\leq T^\vartheta}}\frac{\alpha_{-1}(m)\alpha_{-1}(n)(m,n)}{mn}
\\
&\qquad+O\bigg(T\mathscr{L}^2\sum_{lmn\in S(X)}\frac{d(lm)d(ln)\log mn}{lmn}\bigg)+O_B(T\mathscr{L}^{-B})+O_\varepsilon(T^{2\vartheta+\varepsilon})
\\
&=\frac{T\mathscr{L}^3}{3e^{\gamma_0}\log X}\big(1+O\big((\log X)^{-1}\big)\big)+O_\varepsilon(T^{2\vartheta+\varepsilon}).
\end{split}
\end{equation}

To obtain the estimates for \eqref{206} and \eqref{207} from \eqref{208} and \eqref{209}, we use 
the well known approximation 
\begin{equation} \label{stirling_chi}
\frac{\chi'}{\chi}(\tfrac{1}{2}+it)=-\log\frac{t}{2\pi}+O(t^{-1})\quad(\text{for }t\geq1)
\end{equation}
and  integration by parts. In this way we deduce that
\begin{equation}\label{227}
J_3=-\frac{T\mathscr{L}}{2\pi}\frac{\mathscr{L}^3}{e^{\gamma_0}\log X}\big(1+O\big((\log X)^{-1}\big)\big)+O_\varepsilon(T^{2\vartheta+\varepsilon})
\end{equation}
and 
\begin{equation}\label{228}
J_4=\frac{T\mathscr{L}}{2\pi}\frac{\mathscr{L}^3}{e^{\gamma_0}\log X}\big(1+O\big((\log X)^{-1}\big)\big)+O_\varepsilon(T^{2\vartheta+\varepsilon}).
\end{equation}

\subsection{The evaluation of  $J_1$ and $J_2$}

Note that $J_1+J_2$ equals 
\begin{equation*}
-2 \, \Re \bigg\{\frac{1}{2\pi i}\int_{c+i}^{c+iT}\chi(1-s)\bigg(\frac{\zeta'}{\zeta}(s)\zeta'(s)^2\sum_{\substack{m\in S(X)\\ m\leq T^\vartheta}}\frac{\alpha_{-1}(m)}{m^s}\bigg)\bigg(\sum_{\substack{n\in S(X)\\ n\leq T^\vartheta}}\frac{\alpha_{-1}(n)}{n^{1-s}}\bigg)ds\bigg\}.
\end{equation*}
By Lemma~\ref{lem 2}, we find that 
\begin{equation*}
J_1+J_2=-2 \, \Re \bigg\{\sum_{\substack{n\in S(X)\\ n\leq T^\vartheta}}\frac{\alpha_{-1}(n)}{n}
\sum_{m\leq nT/2\pi}a(m)e ( -m/n )\bigg\}+ O_\varepsilon(T^{1/2+\vartheta+\varepsilon}),
\end{equation*}
where the arithmetic function $a(m)$ is defined by
\begin{equation}\label{211}
\frac{\zeta'}{\zeta}(s)\zeta'(s)^2\sum_{\substack{m\in S(X)\\ m\leq T^\vartheta}}
\frac{\alpha_{-1}(m)}{m^s}=\sum_{m=1}^{\infty}\frac{a(m)}{m^s}
\end{equation}
for $\Re(s)>1$. By the work of Conrey, Ghosh and Gonek (see Sections 5 and  6 and (8.2) of \cite{CGG1}), and of Bui and Heath-Brown \cite{BH-B}, we have 
\begin{equation*}
J_1+J_2 = M_R + E_R+O_\varepsilon(T^{1/2+\vartheta+\varepsilon}),
\end{equation*}
where 
\begin{equation}\label{210}
M_R=-2\sum_{\substack{ln\in S(X)\\ ln\leq T^\vartheta}}\frac{\alpha_{-1}(ln)}{ln}\frac{\mu(n)}{\varphi(n)}\sum_{\substack{m\leq nT/2\pi\\(m,n)=1}}a(lm)
\end{equation}
and
\begin{equation}\label{229}
E_R\ll_{c,B,\varepsilon} T\exp\big(-c\sqrt{\log T}\big)+T\mathscr{L}^{-B}+T^{5/6+\vartheta/3+\varepsilon}
\end{equation}
for some absolute constant $c>0$, and for any $B>0$.

Write
\begin{equation*}
\left(-\frac{\zeta'}{\zeta}(s)\right)^j=\sum_{m=1}^{\infty}\frac{\Lambda_j(m)}{m^s}\quad\textrm{and}\quad-\frac{\zeta'}{\zeta}(s)\zeta'(s)^2=\sum_{m=1}^{\infty}\frac{g(m)}{m^s}
\end{equation*}
for $\Re(s)>1$. From \eqref{211} and Lemma 3.3, we see that
\begin{equation*}
a(lm)=-\sum_{\substack{l=l_1l_2\\m=m_1m_2\\(m_2,l_1)=1}}g(l_1m_1)\alpha_{-1}(l_2m_2),
\end{equation*}
and thus
\begin{equation*}
\begin{split}
\sum_{\substack{m\leq nT/2\pi\\(m,n)=1}}a(lm)=-\sum_{l=l_1l_2} \; \sum_{\substack{l_2m_2\in S(X)\\l_2m_2\leq T^\vartheta\\(m_2,l_1n)=1}}\alpha_{-1}(l_2m_2)\sum_{\substack{m_1\leq nT/2\pi m_2\\(m_1,n)=1}}g(l_1m_1).
\end{split}
\end{equation*}
The innermost sum on the right-hand side has been evaluated by Conrey, Ghosh and Gonek. By Lemma A  of \cite{CGG1}, the sum over $m_1$ is
\begin{equation} \label{CGG_est}
\begin{split}
&=\frac{nT}{2\pi m_2}\frac{\varphi(n)^2}{n^2}\sum_{j=0}^{3}\binom{3}{j}\beta_{j}(l_1)\delta(l_1)\frac{(\log nT/2\pi m_2)^{j+1}}{(j+1)!}+O\big(nT\mathscr{L}^3d(l_1)/m_2\big) \\
&=\frac{T\mathscr{L}^4}{48\pi}\frac{\varphi(n)^2\delta(l_1)}{m_2n}+O\big(T\mathscr{L}^3\varphi(n)d(l_1)(\log l_1n)/m_2\big),
\end{split}
\end{equation}
where $\delta(l)=\prod_{p|l}(2-1/p)$ and $\beta_j(l)=\sum_{d|l}\Lambda_{3-j}(d)/\delta(d)$. 
We insert this estimate into  \eqref{210}. The contribution of the big-$O$ term in the last line of \eqref{CGG_est} to \eqref{210} is
\begin{equation*}
\begin{split}
\ll T\mathscr{L}^3\sum_{l_1l_2mn\in S(X)} \frac{d(l_1l_2n)d(l_2m)d(l_1)\log (l_1n)}{l_1l_2mn}\ll T\mathscr{L}^3\bigg(\sum_{n\in S(X)}\frac{d(n)^2\log n}{n}\bigg)^4.
\end{split}
\end{equation*}
By the same method we used to obtain the estimate in \eqref{212}, the sum over $n$ on the right-hand side is $\ll (\log X)^5$. 
Thus, the contribution from the big-$O$ term is $O(T\mathscr{L}^3(\log X)^{20})$. We therefore have that
\begin{equation*}
\begin{split}
M_R=\frac{T\mathscr{L}^4}{24\pi}\sum_{\substack{l_1l_2n\in S(X)\\ l_1l_2n\leq T^\vartheta}}\sum_{\substack{m\in S(X)\\l_2m\leq T^\vartheta\\(m,l_1n)=1}}\frac{\alpha_{-1}(l_2m)\alpha_{-1}(l_1l_2n)\mu(n)\varphi(n)\delta(l_1)}{l_1l_2mn^2}+O(T\mathscr{L}^3(\log X)^{20}).
\end{split}
\end{equation*}




Next we   show that we may extend the sums to all products $l_1l_2mn\in S(X)$ with $(m,l_1n)=1$ with an acceptable error term. This follows from ``Rankin's trick'', for  we have
\begin{eqnarray*}
\sum_{\substack{l_1l_2mn\in S(X)\\l_1l_2mn>T^\vartheta}}\frac{d(l_2m)d(l_1l_2n)d(l_1)}{l_1l_2mn}&\ll& \sum_{l_1l_2mn\in S(X)}\frac{d(l_2m)d(l_1l_2n)d(l_1)}{l_1l_2mn}\bigg(\frac{l_1l_2mn}{T^\vartheta}\bigg)^{1/4}
\\
&\ll& T^{-\vartheta/4}\bigg(\sum_{n\in S(X)}\frac{d(n)^2}{n^{3/4}}\bigg)^4\ll T^{-\vartheta/4}\prod_{p\leq X}\bigg(1-\frac{1}{p^{3/4}}\bigg)^{-16}
\\
&\ll &T^{-\vartheta/4}e^{100 X^{1/4}/\log X}\ll T^{-\vartheta/5} 
\end{eqnarray*}
since $X=O((\log T)^{2-\varepsilon})$.
Hence,   writing $n$ for $l_1 n$ and $l$ for $l_2$, we have
\begin{equation}\label{M with error term}
M_R=\frac{T\mathscr{L}^4}{24\pi}\sum_{\substack{lmn\in S(X)\\ (m,n)=1}}\frac{\alpha_{-1}(lm)\alpha_{-1}(ln)g(n)}{lmn}+O(T\mathscr{L}^3(\log X)^{20}),
\end{equation}
where
\begin{equation*}
g(n)=\sum_{d|n}\frac{\mu(d)\varphi(d)\delta(n/d)}{d}.
\end{equation*}

Let $P=\prod_{p\leq X}p$. Since $\alpha_{-1}(n)=0$ if $n$ is not a cube-free integer, we can restrict the summation over $l$ to summation over $l=u_{1}u_{2}^{2}$, where $u_{1}|P$, and $u_{2}|(P/u_{1})$. The summation over $m$ and $n$ can also be restricted to $(m,u_{2})=(n,u_{2})=1$, since otherwise $\alpha_{-1}(lm)\alpha_{-1}(ln)=0$. Thus, apart from the big-$O$ term in \eqref{M with error term}, we see that $M_R$ equals
\begin{equation*}
\frac{T\mathscr{L}^4}{24\pi}\sum_{u_{1}|P}\frac{1}{u_{1}}
\sum_{u_{2}|(P/u_{1})}\frac{\alpha_{-1}(u_{2}^{2})^2}{u_{2}^{2}}
\sum_{\substack{m\in S(X)\\(m,u_{2})=1}}\sum_{\substack{n\in S(X)\\(n,u_{2}m)=1}}
\frac{\alpha_{-1}(u_{1}m)\alpha_{-1}(u_{1}n)g(n)}{mn}.
\end{equation*}
Arguing similarly, we see that if $r=(u_{1},m)$ and $m=rm_{1}$, then we can 
assume that $(r,m_{1})=1$ so that $(u_{1},m_{1})=1$. Consequently, the summation over 
$m$ can be replaced by
\begin{equation*}
\sum_{r|u_{1}}\sum_{\substack{m_{1}\in S(X)\\(m_{1},u_{1}u_{2})=1}}.
\end{equation*}
Similarly, for $s=(u_{1},n)$ and $n=sn_{1}$, we can sum over $(u_{1},n_{1})=1$. 
The condition $(m,n)=1$ is equivalent to $(m_{1},n_{1})=(m_{1},s)=(r,n_{1})=(r,s)=1$. 
Now, $(r,s)=1$ if and only if $s|(u_{1}/r)$. Also, $(m_{1},s)=1$ and $(n_1,r)=1$ 
are implied by $(m_{1}n_{1},u_{1})=1$. Thus, $M_R$ equals
\begin{equation*}
\begin{split}
& \frac{T\mathscr{L}^4}{24\pi}\sum_{u_{1}|P}\frac{1}{u_{1}}\sum_{u_{2}|(P/u_{1})}
\frac{\alpha_{-1}(u_{2}^{2})^2}{u_{2}^{2}}\sum_{r|u_{1}}\sum_{\substack{m_{1}\in S(X)\\(m_{1},u_{1}u_{2})=1}}\\
&\qquad\qquad\qquad\qquad\sum_{s|(u_{1}/r)}\sum_{\substack{n_{1}\in S(X)\\(n_{1},u_{1}u_{2}m_{1})=1}}\frac{\alpha_{-1}(u_{1}rm_{1})\alpha_{-1}(u_{1}sn_{1})g(sn_{1})}{rsm_{1}n_{1}}\nonumber\\
&=\frac{T\mathscr{L}^4}{24\pi}\sum_{u_{1}|P}\frac{\alpha_{-1}(u_{1})^2}{u_{1}}\sum_{u_{2}|(P/u_{1})}\frac{\alpha_{-1}(u_{2}^{2})^2}{u_{2}^{2}}\sum_{r|u_{1}}\frac{\alpha_{-1}(r^2)}{\alpha_{-1}(r)r}\sum_{s|(u_{1}/r)}\frac{\alpha_{-1}(s^2)g(s)}{\alpha_{-1}(s)s}\nonumber\\
&\qquad\qquad\qquad\qquad\sum_{\substack{m_{1}\in S(X)\\(m_{1},u_{1}u_{2})=1}}\frac{\alpha_{-1}(m_{1})}{m_{1}}\sum_{\substack{n_{1}\in S(X)\\(n_{1},u_{1}u_{2}m_{1})=1}}\frac{\alpha_{-1}(n_{1})g(n_{1})}{n_{1}}.
\end{split}
\end{equation*}
Since $m_1$ and $n_1$ make no contribution unless they are cube-free,  this last expression is equal to
\begin{equation}\label{226}
\begin{split}
&\frac{T\mathscr{L}^4}{24\pi}\sum_{u_{1}|P}\frac{\alpha_{-1}(u_{1})^2}{u_{1}}\sum_{u_{2}|(P/u_{1})}\frac{\alpha_{-1}(u_{2}^{2})^2}{u_{2}^{2}}\sum_{r|u_{1}}\frac{\alpha_{-1}(r^2)}{\alpha_{-1}(r)r}\sum_{s|(u_{1}/r)}\frac{\alpha_{-1}(s^2)g(s)}{\alpha_{-1}(s)s}
\\
&\qquad\qquad\qquad\sum_{m_{1}|(P/u_{1}u_{2})^2}\frac{\alpha_{-1}(m_{1})}{m_{1}}\sum_{n_{1}|(P/u_{1}u_{2}m_{1})^{2}}\frac{\alpha_{-1}(n_{1})g(n_{1})}{n_{1}}.
\end{split}
\end{equation}

Next we define the following multiplicative functions:
\begin{equation*}
\begin{split}
&T_{1}(n)=\sum_{d|n}\frac{\alpha_{-1}(d)g(d)}{d},\quad T_{2}(n)=\sum_{d|n}\frac{\alpha_{-1}(d)}{dT_{1}(d^2)},\\
&T_{3}(n)=\sum_{d|n}\frac{\alpha_{-1}(d^2)g(d)}{\alpha_{-1}(d)d},\quad T_{4}(n)=\sum_{d|n}\frac{\alpha_{-1}(d^2)}{\alpha_{-1}(d)dT_{3}(d)},\\
&T_{5}(n)=\sum_{d|n}\frac{\alpha_{-1}(d^2)^{2}}{d^{2}T_{1}(d^2)T_{2}(d^2)} \quad \text{and} \quad T_{6}(n)=\sum_{d|n}\frac{\alpha_{-1}(d)^{2}T_{3}(d)T_{4}(d)}{dT_{1}(d^2)T_{2}(d^2)T_{5}(d)}.
\end{split}
\end{equation*}
The  sum over $n_1$ in \eqref{226} equals
\begin{equation*}
T_{1}\big((P/u_{1}u_{2}m_{1})^2\big)=\frac{T_{1}(P^2)}{T_{1}(u_{1}^{2})T_{1}(u_{2}^{2})T_{1}(m_{1}^{2})},
\end{equation*}
and therefore the double summation over $m_{1}$ and $n_{1}$ in \eqref{226} is equal to 
\begin{equation*}
\frac{T_{1}(P^2)}{T_{1}(u_{1}^{2})T_{1}(u_{2}^{2})}T_{2}\big((P/u_{1}u_{2})^2\big)=\frac{T_{1}(P^2)T_{2}(P^2)}{T_{1}(u_{1}^{2})T_{2}(u_{1}^{2})T_{1}(u_{2}^{2})T_{2}(u_{2}^{2})}.
\end{equation*}
Similarly, the summation over $r$ and $s$ in \eqref{226} is
\begin{equation*}
T_3(u_1)T_4(u_1).
\end{equation*}
It follows that 
\begin{eqnarray*}
M_R&=&\frac{T\mathscr{L}^4}{24\pi}T_{1}(P^2)T_{2}(P^2)\sum_{u_{1}|P}\frac{\alpha_{-1}(u_{1})^2T_{3}(u_{1})T_{4}(u_{1})}{u_{1}T_{1}(u_{1}^{2})T_{2}(u_{1}^{2})}\sum_{u_{2}|(P/u_{1})}\frac{\alpha_{-1}(u_{2}^{2})^2}{u_{2}^{2}T_{1}(u_{1}^{2})T_{2}(u_{1}^{2})}\nonumber\\
&=&\frac{T\mathscr{L}^4}{24\pi}T_{1}(P^2)T_{2}(P^2)T_{5}(P)T_{6}(P)\nonumber\\
&=&\frac{T\mathscr{L}^4}{24\pi}\prod_{p\leq X}\bigg(T_{1}(p^2)T_{2}(p^2)T_{5}(p)+\frac{\alpha_{-1}(p)^{2}T_{3}(p)T_{4}(p)}{p}\bigg).
\end{eqnarray*}
To simplify this expression, first note that 
\begin{equation*}
g(p)=1\quad\textrm{ and }\quad g(p^2)=\frac{2}{p}-\frac{1}{p^2}.
\end{equation*}
Moreover, $\alpha_{-1}(p)=-1$ for all $p\leq X$, so  
\begin{equation*}
\frac{\alpha_{-1}(p)^{2}T_{3}(p)T_{4}(p)}{p}=\frac{T_{3}(p)}{p}-\frac{\alpha_{-1}(p^2)}{p^2}=\frac{1}{p}-\frac{2\alpha_{-1}(p^2)}{p^2},
\end{equation*}
and
\begin{equation*}
\begin{split}
T_{1}(p^2)T_{2}(p^2)T_{5}(p)&=T_{1}(p^2)T_{2}(p^2)+\frac{\alpha_{-1}(p^2)^{2}}{p^2}
\\
&=T_{1}(p^2)-\frac{1}{p}+\frac{\alpha_{-1}(p^2)}{p^2}+\frac{\alpha_{-1}(p^2)^{2}}{p^2}
\\
&=1-\frac{2}{p}+\frac{\alpha_{-1}(p^2)\big(1+g(p^2)+\alpha_{-1}(p^2)\big)}{p^2}.
\end{split}
\end{equation*}
Since we also have that $\alpha_{-1}(p^2)=0$ if $p\leq\sqrt{X}$, we see that
\begin{equation*}
\begin{split}
T_{1}(p^2)T_{2}(p^2)&T_{5}(p)+\frac{\alpha_{-1}(p)^{2}T_{3}(p)T_{4}(p)}{p}
\\
&=\left\{ \begin{array}{ll}
1-1/p, \qquad\qquad\quad \qquad\textrm{if $p\leq\sqrt{X}$,}\\
1-1/p+O(1/p^2), \ \qquad\textrm{if $\sqrt{X}<p\leq X$.}
\end{array} \right.
\end{split}
\end{equation*}
Collecting these estimates, we now have, apart from the big-$O$ term in \eqref{M with error term}, that
\begin{equation*}
\begin{split}
M_R&=\frac{T\mathscr{L}^4}{24\pi}\prod_{p\leq\sqrt{X}}\bigg(1-\frac{1}{p}\bigg)\prod_{\sqrt{X}<p\leq X}\bigg(1-\frac{1}{p}+O(1/p^2)\bigg)
\\
&=\frac{T\mathscr{L}}{2\pi}\frac{\mathscr{L}^3}{12e^{\gamma_0}\log X}\big(1+O\big((\log X)^{-1}\big)\big).
\end{split}
\end{equation*}
Combining this expression with \eqref{201}, \eqref{J sum}, \eqref{202}, \eqref{227}, \eqref{228}, \eqref{229}, and \eqref{M with error term}, 
we obtain
\begin{equation*}
\begin{split}
\frac{1}{N(T)}\sum_{0<\gamma\leq T}\big|\zeta'(\rho)P_X(\rho)^{-1}\big|^2&=\frac{\mathscr{L}^3}{12e^{\gamma_0}\log X}\big(1+O\big((\log X)^{-1}\big)\big)
\\
&\quad +O_c\big(\exp\big(-c\sqrt{\log T}\big)\big) +O_B(\mathscr{L}^{-B}) 
\\
&\quad +O_\varepsilon\big(T^{-1/2+\vartheta+\varepsilon}+T^{-1+2\vartheta+\varepsilon}+T^{-1/6+\vartheta/3+\varepsilon}\big).
\end{split}
\end{equation*}
Theorem \ref{2nd_moment} now follows by choosing any $\vartheta<1/2$.

\section{The twisted moment conjectures}

In this section, we use a modification of the recipe in \cite{CFKRS, CS} to formulate a conjecture for the discrete moments of $Z_{X}'(\rho)$. We start by considering the twisted $2k$-th moment of the derivative of the Riemann zeta-function, that is
\begin{displaymath}
I_{2k}(m,n)=\sum_{0<\gamma\leq T}\big|\zeta'(\rho)\big|^{2k}\bigg(\frac{m}{n}\bigg)^{-i\gamma}.
\end{displaymath}
We assume RH and, for simplicity, we assume that $(m,n)=1$. Using Cauchy's theorem, we may write this sum as a contour integral; namely
\begin{displaymath}
I_{2k}(m,n)=\frac{1}{2\pi i}\int_{\mathscr{C}}\frac{\zeta'(s)}{\zeta(s)}\zeta'(s)^{k}\zeta'(1-s)^{k}\bigg(\frac{m}{n}\bigg)^{-s+1/2}ds,
\end{displaymath}
with the contour $\mathscr{C}$ running from $1-c+i$ to $c+i$, $c+iT$ and $1-c+iT$, where as before $c=1+\mathscr{L}^{-1}$. Using standard estimates for the integrand, we can show that the contribution from the horizontal segments of the contour is negligible. Therefore, it suffices to estimate the right-hand and left-hand portions of the contour, $I_{2k,R}(m,n)$ and $I_{2k,L}(m,n)$, say. We first examine the integral from $c+i$ to $c+iT$, which is
\begin{eqnarray*}
I_{2k,R}(m,n)&=&\frac{1}{2\pi}\int_{1}^{T}\frac{\zeta'(c+it)}{\zeta(c+it)}\zeta'(c+it)^{k}\zeta'(1-c-it)^{k}\bigg(\frac{m}{n}\bigg)^{-c-it+1/2}dt\\
&=&\frac{d}{d\alpha_1}\ldots\frac{d}{d\alpha_{k+1}}\frac{d}{d\beta_1}\ldots\frac{d}{d\beta_{k}}\frac{1}{2\pi}\int_{1}^{T}\frac{\zeta(c+it+\alpha_{k+1})}{\zeta(c+it)}\\
&&\quad \quad\times \prod_{j=1}^{k}\bigg(\zeta(c+it+\alpha_j)\zeta(1-c-it+\beta_j)\bigg)\bigg(\frac{m}{n}\bigg)^{-c-it+1/2}dt\Bigg|_{\underline{\alpha}=\underline{\beta}=0}.
\end{eqnarray*}
Following the recipe outlined in \cite{CFKRS, CS}, we replace each of the zeta-functions in the numerator by 
\begin{displaymath}
\zeta(s)\sim\sum_{n\leq\sqrt{t/2\pi}}\frac{1}{n^s}+\chi(s)\sum_{n\leq\sqrt{t/2\pi}}\frac{1}{n^{1-s}},
\end{displaymath}
and we replace the zeta-function in the denominator by
\begin{displaymath}
\frac{1}{\zeta(s)}=\sum_{n=1}^{\infty}\frac{\mu(n)}{n^s}.
\end{displaymath}
Multiplying out the various sums, we obtain $2^{2k+1}$ terms in the integrand. We note that Stirling's formula for the Gamma function implies that
\begin{equation}\label{5}
\chi(s+\alpha)\chi(1-s+\beta)=\bigg(\frac{t}{2\pi}\bigg)^{-\alpha-\beta}\bigg(1+O\bigg(\frac{1}{|t|\!+\!1}\bigg)\bigg)
\end{equation}
as $t \to \infty.$ We only keep the terms with the same number of $\chi$ factors coming from $\zeta(s)$ and  from $\zeta(1-s)$. Consider  the term coming from the product of  the first term  of each approximate functional equation, namely
\begin{displaymath}
\sum_{\substack{a_1,\ldots,a_{k+2}\\b_1,\ldots,b_k}}\frac{\mu(a_{k+2})}{a_1^{\alpha_1}\ldots a_{k+1}^{\alpha_{k+1}}b_1^{1+\beta_1}\ldots b_{k}^{1+\beta_{k}}}\bigg(\frac{a_1\ldots a_{k+2}}{b_1\ldots b_k}\bigg)^{-c-it}\bigg(\frac{m}{n}\bigg)^{-c-it+1/2}.
\end{displaymath}
Averaging over $t$, only the diagonal terms $a_1\ldots a_{k+2}m=b_1\ldots b_kn$ are retained and we obtain
\begin{equation}\label{1}
\int_{1}^{T}\sum_{am=bn}\frac{A_{\underline{\alpha}}(a)B_{\underline{\beta}}(b)}{\sqrt{ab}}dt,
\end{equation}
where
\begin{displaymath}
A_{\underline{\alpha}}(a)=\sum_{a_1\ldots a_{k+2}=a}\frac{\mu(a_{k+2})}{a_1^{\alpha_1}\ldots a_{k+1}^{\alpha_{k+1}}},
\end{displaymath}
and
\begin{displaymath}
B_{\underline{\beta}}(b)=\sum_{b_1\ldots b_{k}=b}\frac{1}{b_1^{\beta_1}\ldots b_{k}^{\beta_{k}}}.
\end{displaymath}
Since $(m,n)=1$, the only solutions of $am=bn$ are $a=un$ and $b=um$. Thus, since $A_{\underline{\alpha}}(a)$ and $B_{\underline{\beta}}(b)$ are multiplicative functions, the integral in $\eqref{1}$ equals
\begin{eqnarray*}
&& \frac{1}{\sqrt{mn}}\int_{1}^{T}\sum_{u=1}^{\infty}
 \frac{A_{\underline{\alpha}}(un)B_{\underline{\beta}}(um)}{u}dt=\frac{1}{\sqrt{mn}}\int_{1}^{T}\sum_{u=1}^{\infty}
 \frac{A_{\underline{\alpha}}(u)B_{\underline{\beta}}(u)}{u} \\
&&\qquad\qquad\times
  \prod_{\substack{p^{m_p}||m\\p^{n_p}||n}}\bigg(\frac{\sum_{j=0}^{\infty}A_{\underline{\alpha}}(p^{j+n_p})B_{\underline{\beta}}(p^j)/p^j}{\sum_{j=0}^{\infty}A_{\underline{\alpha}}(p^{j})
  B_{\underline{\beta}}(p^j)/p^j}\frac{\sum_{j=0}^{\infty}A_{\underline{\alpha}}(p^{j})B_{\underline{\beta}}(p^{j+m_p})/p^j}{\sum_{j=0}^{\infty}A_{\underline{\alpha}}(p^{j})B_{\underline{\beta}}(p^j)/p^j}\bigg)dt.
\end{eqnarray*}
We denote the integrand on the right-hand side of the above equation by $T_{\underline{\alpha},\underline{\beta}}(m,n)$, and we denote the product over primes in this integrand by $C_{\underline{\alpha},\underline{\beta}}(m,n)$. Now the sum over $u$ in $T_{\underline{\alpha},\underline{\beta}}(m,n)$  is
\begin{equation*}
\begin{split}
\sum_{u=1}^{\infty}\frac{A_{\underline{\alpha}}(u)B_{\underline{\beta}}(u)}{u} =\prod_{p}\bigg(\sum_{\sum_{j=1}^{k+2}a_j=\sum_{j=1}^{k}b_j}\frac{\mu(p^{a_{k+2}})}{p^{\sum_{j=1}^{k+1}(1/2+\alpha_j)a_j+a_{k+2}/2+\sum_{j=1}^{k}(1/2+\beta_j)b_j}}\bigg).
\end{split}
\end{equation*}
Taking out the divergent terms from the above formula in the form of zeta-functions, the integrand $T_{\underline{\alpha},\underline{\beta}}(m,n)$ equals
\begin{eqnarray*}
&&\frac{\prod_{\substack{1\leq i\leq k+1\\1\leq j\leq k}}\zeta(1+\alpha_{i}+\beta_j)}{\prod_{1\leq j\leq k}\zeta(1+\beta_j)} \prod_{p}
\Bigg(\prod_{\substack{1\leq i\leq k+1\\1\leq j\leq k}}\bigg(1-\frac{1}{p^{1+\alpha_{i}+\beta_j}}\bigg)\prod_{1\leq j\leq k}\bigg(1-\frac{1}{p^{1+\beta_j}}\bigg)^{-1}\\
&&\qquad\qquad\times\sum_{\sum_{j=1}^{k+2}a_j=\sum_{j=1}^{k}b_j}\frac{\mu(p^{a_{k+2}})}{p^{\sum_{j=1}^{k+1}(1/2+\alpha_j)a_j+a_{k+2}/2+\sum_{j=1}^{k}(1/2+\beta_j)b_j}}\Bigg)C_{\underline{\alpha},\underline{\beta}}(m,n).
\end{eqnarray*}
We handle the other terms which arise from multiplying out the approximate functional equations in a similar manner, but we also take into account the asymptotic formula  \eqref{5}. Adding the resulting terms, we obtain
\begin{equation*}
\begin{split}
I_{2k,R}(m,n)&=\frac{d}{d\alpha_1}\ldots\frac{d}{d\alpha_{k+1}}\frac{d}{d\beta_1}\ldots\frac{d}{d\beta_{k}}\frac{1}{2\pi\sqrt{mn}}\\
& \quad \times \int_{1}^{T}\sum_{0\leq j\leq k}\sum_{\substack{P\subset\{\alpha_1,\ldots,\alpha_{k+1}\}\\Q\subset\{\beta_1,\ldots,\beta_k\}\\|P|=|Q|=j}}T_{\underline{\alpha}_P,\underline{\beta}_Q}(m,n)\bigg(\frac{t}{2\pi}\bigg)^{-P-Q}dt\bigg|_{\underline{\alpha}=\underline{\beta}=0}+O_{k,\varepsilon}(T^{1/2+\varepsilon}),
\end{split}
\end{equation*}
where if $P=\{\alpha_{u_1},\ldots,\alpha_{u_j}\}$ and $Q=\{\beta_{v_1},\ldots,\beta_{v_j}\}$ with $u_1<\ldots<u_j$ and $v_1<\ldots<v_j$, then $(\underline{\alpha}_P,\underline{\beta}_Q)$ is the $(2k\!+\!1)$-tuple obtained from 
\[
(\alpha_1,\ldots,\alpha_{k+1},\beta_1,\ldots,\beta_k)
\] 
by replacing $\alpha_{u_r}$ with $-\beta_{v_r}$ and replacing $\beta_{v_r}$ with $-\alpha_{u_r}$ for all $1\leq r\leq j$. Here $(t/2\pi)^{-P-Q}$ stands for 
\[
(t/2\pi)^{-\sum_{x\in P}x-\sum_{y\in Q}y}.
\] 
There is a concise way to write these $\binom{2k+1}{k}$ terms as a contour integral (see \cite{CFKRS}), namely $I_{2k,R}(m,n)$ equals
\begin{eqnarray*}
&&\frac{d}{d\alpha_1}\ldots\frac{d}{d\alpha_{k+1}}\frac{d}{d\beta_1}\ldots\frac{d}{d\beta_{k}}\frac{1}{2\pi\sqrt{mn}}\int_{1}^{T}\bigg(\frac{t}{2\pi}\bigg)^{\frac{-\sum_{j}\alpha_j-\sum_{j}\beta_j}{2}}\frac{1}{(k+1)!k!(2\pi i)^{2k+1}}\\
 &&\quad\quad\times \oint\ldots\oint\bigg(\frac{t}{2\pi}\bigg)^{\frac{\sum_{j}s_j-\sum_{j}z_j}{2}}\frac{T_{\underline{s},-\underline{z}}(m,n)\Delta(s_1,\ldots,s_{k+1},z_1,\ldots,z_k)^2}{\prod_{i,j}(s_i-\alpha_j)\prod_{i,j}(s_i+\beta_j)\prod_{i,j}(z_i-\alpha_j)\prod_{i,j}(z_i+\beta_j)}\\
 &&\qquad\qquad\qquad \times ds_1\ldots ds_{k+1}dz_1\ldots dz_kdt\Big|_{\underline{\alpha}=\underline{\beta}=0}+O_{k,\varepsilon}(T^{1/2+\varepsilon}),
\end{eqnarray*}
where $\Delta(.)$ is the Vandermonde function and the paths of integration are small circles around the poles $\alpha_j$ and $-\beta_j$. We observe that
\begin{equation}\label{6}
\frac{d}{d\alpha}\frac{e^{-a\alpha}}{\prod_{j=1}^{n}(z_j-\alpha)}\bigg|_{\alpha=0}=\frac{1}{\prod_{j=1}^{n}z_j}\bigg(\sum_{j=1}^{n}\frac{1}{z_j}-a\bigg)
\end{equation}
and
\begin{equation}\label{7}
\frac{d}{d\beta}\frac{e^{-a\beta}}{\prod_{j=1}^{n}(z_j+\beta)}\bigg|_{\beta=0}=\frac{1}{\prod_{j=1}^{n}z_j}\bigg(-\sum_{j=1}^{n}\frac{1}{z_j}-a\bigg).
\end{equation}
Thus
\begin{eqnarray*}
&&I_{2k,R}(m,n) =  \frac{1}{2\pi\sqrt{mn}(k+1)!k!(2\pi i)^{2k+1}}\int_{1}^{T}\oint\ldots\oint\bigg(\frac{t}{2\pi}\bigg)^{\frac{\sum_{j}s_j-\sum_{j}z_j}{2}}\nonumber\\
&&\qquad \qquad \times\frac{T_{\underline{s},-\underline{z}}(m,n)\Delta(s_1,\ldots,s_{k+1},z_1,\ldots,z_k)^2}{(\prod_{j=1}^{k+1}s_j\prod_{j=1}^{k}z_j)^{2k+1}}\bigg(-\frac{\mathscr{L}}{2}+\sum_{j=1}^{k+1}\frac{1}{s_j}+\sum_{j=1}^{k}\frac{1}{z_j}\bigg)^{k+1}\nonumber\\
&&\qquad \qquad \times\bigg(-\frac{\mathscr{L}}{2}-\sum_{j=1}^{k+1}\frac{1}{s_j}-\sum_{j=1}^{k}\frac{1}{z_j}\bigg)^{k}ds_1\ldots ds_{k+1}dz_1\ldots dz_kdt+O_{k,\varepsilon}(T^{1/2+\varepsilon}).
\end{eqnarray*}

The contribution from the left-hand side of the contour of integration is
\begin{displaymath}
I_{2k,L}(m,n)=-\frac{1}{2\pi}\int_{1}^{T}\frac{\zeta'(1-c+it)}{\zeta(1-c+it)}\zeta'(1-c+it)^{k}\zeta'(c-it)^{k}\bigg(\frac{m}{n}\bigg)^{c-it-1/2}dt.
\end{displaymath}
By the functional equation for $\zeta'(s)/\zeta(s)$ in \eqref{log_deriv}, we have 
\begin{displaymath}
\frac{\zeta'(1-c+it)}{\zeta(1-c+it)}=\frac{\chi'(1-c+it)}{\chi(1-c+it)}-\frac{\zeta'(c-it)}{\zeta(c-it)}.
\end{displaymath}
Thus,
\begin{equation*}
\begin{split}
I_{2k,L}(m,n)&=-\frac{1}{2\pi i}\int_{1-c+i}^{1-c+iT}\frac{\chi'(s)}{\chi(s)}\zeta'(s)^{k}\zeta'(1-s)^{k}\bigg(\frac{m}{n}\bigg)^{-s+1/2}ds\\
&\quad\quad+\frac{1}{2\pi}\int_{1}^{T}\frac{\zeta'(c-it)}{\zeta(c-it)}\zeta'(c-it)^{k}\zeta'(1-c+it)^{k}\bigg(\frac{m}{n}\bigg)^{c-it-1/2}dt.
\end{split}
\end{equation*}
We note that the second term on the right-hand side is equal to $\overline{I_{2k,R}(n,m)}$. To handle the first term, we may first shift the line of integration to the $\frac{1}{2}$-line with a negligible error. Then, using the approximation for $\chi'(s)/\chi(s)$ in \eqref{stirling_chi}, we   find that this term is roughly equal to
\begin{eqnarray*}
&&\frac{\mathscr{L}}{2\pi}\int_{1}^{T}\zeta'(\tfrac{1}{2}\!+\!it)^{k}\zeta'(\tfrac{1}{2}\!-\!it)^{k}\bigg(\frac{m}{n}\bigg)^{-it}dt=\frac{d}{d\alpha_1}\ldots\frac{d}{d\alpha_{k}}\frac{d}{d\beta_1}\ldots\frac{d}{d\beta_{k}}
  \frac{\mathscr{L}}{2\pi}
  \\
  &&\qquad\qquad\qquad\qquad\quad \quad \times \int_{1}^{T}\prod_{j=1}^{k}\bigg(\zeta(\tfrac{1}{2}\!+\!it\!+\!\alpha_j)\zeta(\tfrac{1}{2}\!-\!it\!+\!\beta_j)\bigg)\bigg(\frac{m}{n}\bigg)^{-it}dt\Bigg|_{\underline{\alpha}=\underline{\beta}=0}.
\end{eqnarray*}
Hughes and Young~\cite{HY} have conjectured that this integral equals
\begin{displaymath}
\frac{1}{\sqrt{mn}}\int_{1}^{T}\Bigg(\sum_{0\leq j\leq k}\sum_{\substack{P\subset\{\alpha_1,\ldots,\alpha_{k}\}\\Q\subset\{\beta_1,\ldots,\beta_k\}\\|P|=|Q|=j}}S_{\underline{\alpha}_P,\underline{\beta}_Q}(m,n)\bigg(\frac{t}{2\pi}\bigg)^{-P-Q}\Bigg)dt+O_{k,\varepsilon}(T^{1/2+\varepsilon}),
\end{displaymath}
where
\begin{equation*}
\begin{split}
&S_{\underline{\alpha}_P,\underline{\beta}_Q}(m,n)=\prod_{1\leq i,j\leq k}\zeta(1+\alpha_{i}+\beta_j)\\
&\qquad\times\prod_{p}\Bigg(\prod_{1\leq i,j\leq k}\bigg(1-\frac{1}{p^{1+\alpha_{i}+\beta_j}}\bigg)\sum_{\sum_{j=1}^{k}a_j=\sum_{j=1}^{k}b_j}\frac{1}{p^{\sum_{j=1}^{k}(1/2+\alpha_j)a_j+(1/2+\beta_j)b_j}}\Bigg)D_{\underline{\alpha},\underline{\beta}}(m,n),
\end{split}
\end{equation*}
with
\begin{equation*}
\begin{split}
D_{\underline{\alpha},\underline{\beta}}(m,n)=\prod_{\substack{p^{m_p}||m\\p^{n_p}||n}}\bigg(\frac{\sum_{j=0}^{\infty}B_{\underline{\alpha}}(p^{j+m_p})B_{\underline{\beta}}(p^j)/p^j}{\sum_{j=0}^{\infty}B_{\underline{\alpha}}(p^{j})B_{\underline{\beta}}(p^j)/p^j}\times\frac{\sum_{j=0}^{\infty}B_{\underline{\alpha}(p^{j})B_{\underline{\beta}}}(p^{j+n_p})/p^j}{\sum_{j=0}^{\infty}B_{\underline{\alpha}}(p^{j})B_{\underline{\beta}}(p^j)/p^j}\bigg).
\end{split}
\end{equation*}
This expression can   be treated as before, that is, by expressing it as a contour integral, 
and using \eqref{6} and \eqref{7}. In this way, we obtain the following conjecture.

\begin{conjecture}\label{con 5} Suppose $m,n \in \mathbb{N}$ with $(m,n)=1$, and $mn \ll_\varepsilon T^{1/2-\varepsilon}$. Then we have
\begin{eqnarray*}
I_{2k}(m,n)
&=&\frac{1}{2\pi\sqrt{mn}(k+1)!k!(2\pi i)^{2k+1}}\int_{1}^{T}\oint\ldots\oint\bigg(\frac{t}{2\pi}\bigg)^{\frac{\sum_{j}s_j-\sum_{j}z_j}{2}}\nonumber\\
&& \qquad \times \frac{\big(T_{\underline{s},-\underline{z}}(m,n)+T_{\underline{s},-\underline{z}}(n,m)\big)\Delta(s_1,\ldots,s_{k+1},z_1,\ldots,z_k)^2}{(\prod_{j=1}^{k+1}s_j\prod_{j=1}^{k}z_j)^{2k+1}}
\\
&&\qquad\times\bigg(-\frac{\mathscr{L}}{2}+\sum_{j=1}^{k+1}\frac{1}{s_j}+\sum_{j=1}^{k}\frac{1}{z_j}\bigg)^{k+1}\\
&&\qquad\times\bigg(-\frac{\mathscr{L}}{2}-\sum_{j=1}^{k+1}\frac{1}{s_j}-\sum_{j=1}^{k}\frac{1}{z_j}\bigg)^{k}ds_1\ldots ds_{k+1}dz_1\ldots dz_kdt\\
&&\ +\frac{\mathscr{L}}{2\pi\sqrt{mn}(k!)^2(2\pi i)^{2k}}\int_{1}^{T}\oint\ldots\oint\bigg(\frac{t}{2\pi}\bigg)^{\frac{\sum_{j}s_j-\sum_{j}z_j}{2}}\\
&& \qquad \times \frac{S_{\underline{s},-\underline{z}}(m,n)\Delta(s_1,\ldots,s_{k},z_1,\ldots,z_k)^2}{(\prod_{j=1}^{k}s_jz_j)^{2k}}
\\
&&\qquad\times\bigg(-\frac{\mathscr{L}}{2}+\sum_{j=1}^{k}\bigg(\frac{1}{s_j}+\frac{1}{z_j}\bigg)\bigg)^{k}\\
&&\qquad\times\bigg(-\frac{\mathscr{L}}{2}-\sum_{j=1}^{k}\bigg(\frac{1}{s_j}
+\frac{1}{z_j}\bigg)\bigg)^{k}ds_1\ldots ds_{k}dz_1\ldots dz_kdt+O_{k,\varepsilon}(T^{1/2+\varepsilon}).
\end{eqnarray*}
\end{conjecture}

We now use Conjecture~\ref{con 5} to give another heuristic argument for Conjectures~\ref{con 3} and \ref{con 4}. Since high moments have much more complicated arithmetic contributions, we shall only treat the case $k=2$. Conjecture~\ref{con 5} asserts that $I_{4}(m,n)$ is asymptotic to $T\mathcal{P}(\mathscr{L})/\sqrt{mn}$, where $\mathcal{P}(x)$ is a polynomial of degree $9$ with coefficients depending on $m$ and $n$. We wish to extract the leading term from this expression. To do this we compute the residues at $s_1=s_{2}=s_3=z_1=z_{2}=0$ of the contour integrals. In this way, we find that
\begin{equation}\label{21}
I_{4}(m,n)=\frac{T\mathscr{L}}{2\pi}\frac{\mathscr{L}^8}{8640\zeta(2)}\frac{\delta(m)\delta(n)}{\sqrt{mn}} +O\big((mn)^{-1/2}d(m)d(n)T\mathscr{L}^8\big),
\end{equation}  
where
\begin{displaymath}
\delta(n)=\prod_{p^{n_p}||n}\bigg(1+n_p\frac{1-1/p}{1+1/p}\bigg).
\end{displaymath}

Using the expression in \eqref{22} with $k=-2,$ we have
\begin{equation}
\begin{split}
\sum_{0<\gamma\leq T}\big|\zeta'(\rho)P_{X}(\rho)^{-1}\big|^4 = \sum_{\substack{mn\in S(X)\\m,n\leq T^{\vartheta}}}\frac{\alpha_{-2}(m)\alpha_{-2}(n)}{\sqrt{mn}}\sum_{0<\gamma\leq T}\big|\zeta'(\rho)\big|^4\bigg(\frac{m}{n}\bigg)^{-i\gamma}+O_{\varepsilon}(T^{1-\varepsilon\vartheta/3}).\label{23}
\end{split}
\end{equation}
It  follows from \eqref{21} that the sum over $m$ and $n$ here equals
\begin{eqnarray}
&&\frac{T\mathscr{L}}{2\pi}\frac{\mathscr{L}^8}{8640\zeta(2)}\sum_{\substack{mn\in S(X)\\m,n\leq T^{\vartheta}}}\frac{\alpha_{-2}(m)\alpha_{-2}(n)\delta\big(m/(m,n)\big)\delta\big(n/(m,n)\big)(m,n)}{mn}\nonumber\\
&&\qquad\qquad+O\bigg(T\mathscr{L}^8\sum_{mn\in S(X)}\frac{d(m)^2d(n)^2(m,n)}{mn}\bigg).\label{524}
\end{eqnarray}
The big-$O$ term is
\begin{displaymath}
\ll T\mathscr{L}^8\sum_{l\in S(X)}\frac{d(l)^4}{l}\bigg(\sum_{m\in S(X)}\frac{d(m)^2}{m}\bigg)^2\ll T\mathscr{L}^8(\log X)^{24},
\end{displaymath}
while the sum over $m$ and $n$ in the main term has been evaluated by Gonek, Hughes and Keating (see pp. 534, 538 of \cite{GHK}) and is
\begin{displaymath}
\sim\frac{\pi^2}{6}(e^{\gamma_0}\log X)^{-4}.
\end{displaymath}
Thus, combining with \eqref{23}, \eqref{524}, and choosing $\vartheta$ sufficiently small, we obtain
\begin{displaymath}
\frac{1}{N(T)}\sum_{0<\gamma\leq T}\big|\zeta'(\rho)P_{X}(\rho)^{-1}\big|^4\sim\frac{1}{8640}\frac{\mathscr{L}^8}{(e^{\gamma_0}\log X)^{4}}.
\end{displaymath}
This heuristic argument provides further evidence for Conjecture~\ref{con 3} and Conjecture \ref{con 4}  in the case $k=2$.

\section*{Acknowledgments}

Work of the second author was supported by National Science Foundation grants DMS-0653809 and DMS-1200582. The third author was supported in part by the NSA Young Investigator Grant H98230-13-1-0217 and an AMS-Simons Travel Grant.

\end{document}